\documentclass{article}

\usepackage[utf8]{inputenc}
\usepackage{graphicx}
\usepackage{epstopdf}
\usepackage{amssymb}
\usepackage{amsfonts}
\usepackage{color}
\usepackage{amsmath}

\def\sn{\hbox{sn}}
\def\sech{\hbox{sech}}

\title{Uniform asymptotic behaviour of Jacobi-$\sn$  near a singular point. The Lost formula from handbooks for elliptic functions}
\author{O.M. Kiselev, ok@ufanet.ru
\\
Institute of Mathematics USC of RAS, Ufa, Russia}

\begin{document}
\maketitle
\begin{abstract}
In this work we construct uniform asymptotic expansion of $\sn(t|m)$ - Jacobi when $m\to1-0$. The constructed expansion is valid over more than a half of period. The turning point is included into the interval of validity for the approximation. In addition we obtain the asymptotic formula for  elliptic integral of the first kind and discuss the differences with the same formula from a handbook.
\end{abstract}
\section{Asymptotic behaviour of elliptic function}

Our goal is to construct  asymptotic formula  for Jacobi elliptic function $\sn(t|m)$, when  $m\to1-0$, which is uniform for all period of the function. The asymptotic expansions for the elliptic functions, when $m\to1-0$, are given in numerous handbooks, see for example  \cite{AbramowitzStegun}. However that expansions not are uniform for very large $t$, when $t=O(T(k))$, where $T(k)$ is a period of the function. An obstacle for the uniformity of the expansions is  the special behaviours of the elliptic functions in neighbourhoods of the turning points $t=T(k)/4+nT/2,\forall n\in \mathbb{Z}$ and far from them.

Let us consider an equation
\begin{equation}
(u')^2=(1-u^2)(1-(1-\epsilon)u^2),\quad 0<\epsilon\ll1.
\label{eqSnJacobi}
\end{equation}
with an initial condition $u(0)=0$. The solution of this Cauchy problem is the Jacobi elliptic function:
$$
u(t,\epsilon)=\sn(t|m),\quad m=1-\epsilon.
$$
The handbook  gives the following approximation (see \cite{AbramowitzStegun}, formula 16.15.1):
$$
\sn(t|1-\epsilon)\sim \tanh(t)+\frac{1}{4}\epsilon\left(\sinh(t)\cosh(t) -t\right)\sech^2(t).
$$

\begin{figure}
\vspace{-7cm}
\includegraphics[scale=0.5]{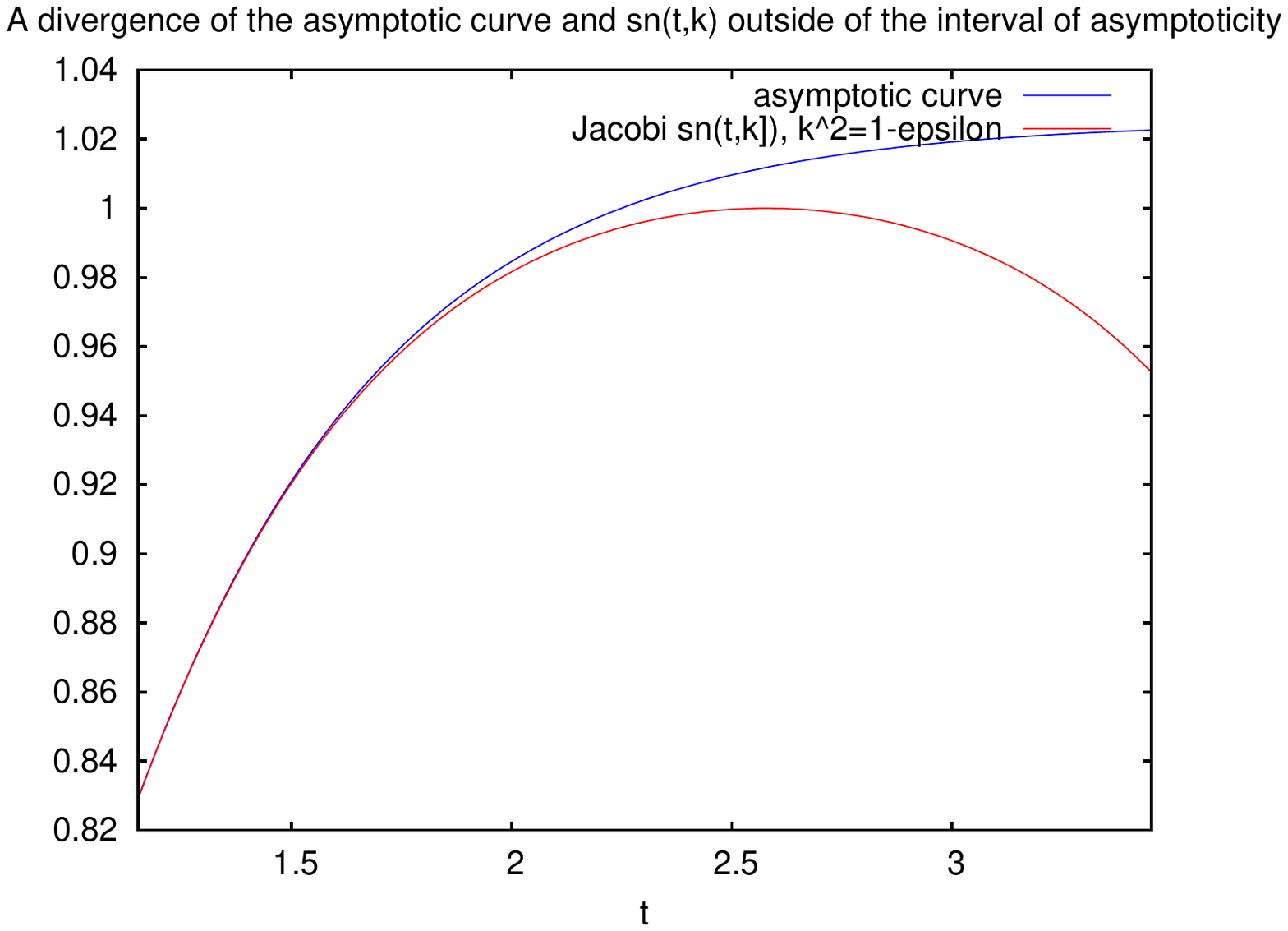}
\caption{The divergence of the asymptotic curve and function $\sn(t|1-\epsilon)$ near the turning point.}
\end{figure}

This approximation is non-periodic, but $\sn(t,1-\epsilon)$ is periodic function with formula for the period:
$$
T(\epsilon)=4\int_0^1\frac{dy}{\sqrt{(1-y^2)(1-(1-\epsilon)y^2)}}\equiv 4 K(1-\epsilon).
$$
The integral in the right-hand side of the formula is the elliptic integral of the  first kind, which is typically denoted by $K(m)$. The handbook \cite{AbramowitzStegun}  gives an polynomial approximation of the integral (formula (17.3.34)):
\begin{eqnarray}
K(m)=
&(1.38662943+ 
0.09666344259 \epsilon+
0.03590092383\epsilon^2+
\nonumber
\\
&
0.03742563713\epsilon^3+
0.01451196212\epsilon^4)+
\nonumber
\\ 
&
(0.5+
0.12498593597\epsilon+
0.06880248576\epsilon^2+
\nonumber
\\
&
0.03328355346\epsilon^3+
0.00441787012
\epsilon^4)\log(1/\epsilon)+e(m),
\nonumber
\\
&
|e(m)|<2\times10^{-8}.
\label{numericApproximationFromAbramowitzStegun}
\end{eqnarray}

In this work we clarify the formula for asymptotic of the elliptic integral of first order and obtain an asymptotic an uniform asymptotic approximation for $\sn(t|1-\epsilon)$. Due to the symmetry
$$
\sn(t|m)=-\sn(-t|m),\quad \sn(t+T/2|m)=-\sn(t|m)
$$
the asymptotic approximation is sufficiently for a half of the period.

\subsection{The asymptotic behaviour of the period}

For small values of $\epsilon$ the elliptic integral can be represented in a form of an integral with weak singularity at $y=1$.

Let us consider the integral as a sum of two integrals over two intervals from zero to small neighbour of $1$ and over small neighbourhood of $y<1$:
\begin{eqnarray*}
\int_0^1\frac{dy}{\sqrt{(1-y^2)(1-(1-\epsilon)y^2)}}=
\int_0^1\frac{dy}{\sqrt{(1+y)(1+\sqrt{1-\epsilon}y)}}\frac{1}{\sqrt{(1-y)(1-\sqrt{1-\epsilon}y)}}.
\end{eqnarray*}
Denote $1-\sqrt{1-\epsilon}=\mu$ then
\begin{eqnarray*}
K(1-\epsilon)=\int_0^1\frac{dy}{\sqrt{(1+y)(1+y-\mu y)}}\frac{1}{\sqrt{(1-y)(1-y+\mu y)}}
\end{eqnarray*}
Now it is convenient to expand the first multiplier into series of $\mu$:
\begin{eqnarray*}
\frac{1}{\sqrt{(1+y)(1+y-\mu y)}}=
&
\frac{1}{y+1}+\frac{\mu y}{2(y+1)^2}+\frac{3\mu^2 y^2}{8(y+1)^3}+
\\
&
\frac{5\mu^3 y^3}{16(y+1)^4}+\frac{35\mu^4 y^4}{128(y+1)^5}+O(\mu^5).
\end{eqnarray*}
Next step is a substitution of this expansion  into integral for $K(1-\epsilon)$. Now the integral should be represented as a sum of integrals. These integrals are obviously calculated by Computer Algebra  System, such as Maxima \cite{maxima}. For example:
\begin{eqnarray*}
I_0=\int_0^1
\frac{dy}{( y+1)\sqrt{(1-y)( 1-( 1-\mu)y)}}=
\\
\frac{\sqrt{4-2\,\mu}\,\log( \mu)}{2\,\mu-4}-\frac{\log( -\mu+2\,\sqrt{4-2\,\mu}+4) \,\sqrt{4-2\,\mu}}{2\,\mu-4}.
\end{eqnarray*}
Similar formulas are obtained for the next integrals:
$$
I_k=\mu^k a_k\int_0^1
\frac{y^k dy}{( y+1)^{k+1}\sqrt{(1-y)( 1-( 1-\mu)y)}}, \quad k=1,1,2,3,4.
$$
Here $a_1=1/2,\,,a_2=3/8,\, a_3=5/16,\,a_4=35/128$. As a result one obtain:
$$
K(1-\epsilon)\sim I_0+I_1+i_2+I_3+I_4,\mu\to0.
$$
This formula in the terms of $\epsilon$ has the following form:
\begin{eqnarray}
K(1-\epsilon)\sim
& 
-\frac{\log(\epsilon)}{2}+2\log(2)+
\left(\frac{-1+2\log(2)}{4}-\frac{\log(\epsilon)}{8}\right) \,\epsilon+
\nonumber
\\
&
\left( \frac{-21+36\,\log( 2) }{128}-\frac{9\,\log( \epsilon)}{128}\right) \,{\epsilon}^{2}+
\nonumber\\
&
\left( \frac{-185+300\,\log( 2) }{1536}-\frac{25\,\log(\epsilon)}{512}\right) \,{\epsilon}^{3}+
\nonumber\\
&
\left( \frac{-18655+29400\,\log( 2) }{196608}-\frac{1225\,\log( \epsilon) }{32768}\right) \,{\epsilon}^{4},\quad \epsilon\to0.
\label{asymptoticOfEllipticIntegralOfFirstKind}
\end{eqnarray}

The similar formula in the numeric form:
\begin{eqnarray}
K(1-\epsilon)\sim
-0.5\,\mathrm{log}\left( \epsilon\right) + 
1.386294361119891+ 
\nonumber
\\
\epsilon\,\left( 0.09657359027997264-0.125\,\mathrm{log}\left( \epsilon\right) \right) +
\nonumber
\\
{\epsilon}^{2}\,\left( 0.03088514453248459-0.0703125\,\mathrm{log}\left( \epsilon\right) \right) 
+
\nonumber
\\
{\epsilon}^{3}\,\left( 0.01493760036978098-0.048828125\,\mathrm{log}\left( \epsilon\right) \right) 
+
\nonumber
\\
{\epsilon}^{4}\,\left( 0.00876631219717606-0.037384033203125\,\mathrm{log}\left( \epsilon\right) \right).
\label{numericOfEllipticIntegralOfFirstKind}
\end{eqnarray}
The differences between (\ref{numericOfEllipticIntegralOfFirstKind}) and (\ref{numericApproximationFromAbramowitzStegun}) can be explained by the different kind of these formulas. The formula (\ref{numericOfEllipticIntegralOfFirstKind}) is asymptotic, but (\ref{numericApproximationFromAbramowitzStegun}) is numerical approximation by polynomials of the value for numerical calculated of the elliptic integral.

\subsection{Asymptotic behaviour on a regular section of trajectory}

Let us construct the solution in the form of expansion on $\epsilon$:
\begin{equation}
u=\tanh(t)+\sum_{n=1}^\infty \epsilon^n u_n(t).
\label{AsymptoticsCloseToSeparatrix}
\end{equation}
Here the primary term of asymptotic expansion is a separatrix of equation (\ref{eqSnJacobi}) as $\epsilon=0$. 

Equations of the high-order terms can be obtained after substituting of  (\ref{AsymptoticsCloseToSeparatrix}) into (\ref{eqSnJacobi}) and collecting  coefficients of  $\epsilon^k$ for $k\i\mathbb{N}$. These equations are combined into a recurrent system. 

For example the equations for  $u_1$ and $u_2$:
\begin{eqnarray*}
\frac{2}{\cosh^2(t)} u_1'+4\frac{\tanh(t)}{\cosh^2(t)}u_1 +\tanh^4(t)-\tanh^2(t)=0
\\
\frac{2}{\cosh^2(t)} u_2'+4\frac{\tanh(t)}{\cosh^2(t)}u_2 +(u_1')^2 +(- 6\tanh^2(t) +2)u_1^2
\\
+(4\tanh^3(t)-2\tanh(t)) u_1=0.
\end{eqnarray*}

The equation for  $n$-th order term has the form:
\begin{eqnarray}
\frac{2}{\cosh^2(t)} u_n'+4\frac{\tanh(t)}{\cosh^2(t)}u_n+\sum_{|\alpha|=n}A_\alpha \tanh^{\alpha_0}(t) u_1^{\alpha_1}\dots u_{n-1}^{\alpha_{n-1}}=0.
\label{EqSnJacobiUn}
\end{eqnarray}
Initial conditions for all corrections are $u_n|_{t=0}=0$.

The equations for the high-order terms have solutions in the form:
$$
u_n=\frac{a_n(t)}{\cosh^2(t)}.
$$
Particularly for $a_1$ one obtains:
$$
a_1'=\sinh^2(t),
$$
It yields:
$$
a_1=\frac{1}{8}\sinh(2t)-\frac{1}{4} t.
$$
The equation for  $a_2$:
$$
a_2'=\frac{1}{64}\cosh(4t)-\frac{5}{32}\cosh(2t)+\frac{1}{8} t \tanh(t)+\frac{1}{16}t^2\sinh^2(t)+\frac{9}{64}.    
$$
It yields:
$$
a_2=-\frac{t^2}{16}\tanh(t)-\frac{1}{256}\sinh(4t)+\frac{5}{64}\sinh(2t)-\frac{9}{64}t.
$$
The high-order terms can be obtained by the same way using  (\ref{EqSnJacobiUn}). In particular for $n$-th order term one obtains:
$$
a_n'=-\frac{1}{2\cosh^{2n-6}(t)}\sum_{|\alpha|=n}A_\alpha \tanh^{\alpha_0}(t)a_1^{\alpha_1}\dots a_{n-1}^{\alpha_{n-1}}
$$

Obvious form of the solutions are very large. Here the asymptotic behaviour as $t\to\pm\infty$ are more important:
$$
a_n=O(e^{\pm 2nt}).
$$
Then the interval of validity for the constructed expansion is:
$$
\epsilon e^{\pm 2t}\ll1,\quad t\ll\mp\frac{1}{2}\log(\epsilon).
$$
As a result we get that the asymptotic expansion of the elliptic functions is valid when $\log(\epsilon)/2\ll t\ll -\log(\epsilon)/2$.

\begin{figure}
\vspace{-10cm}
\includegraphics[scale=0.5]{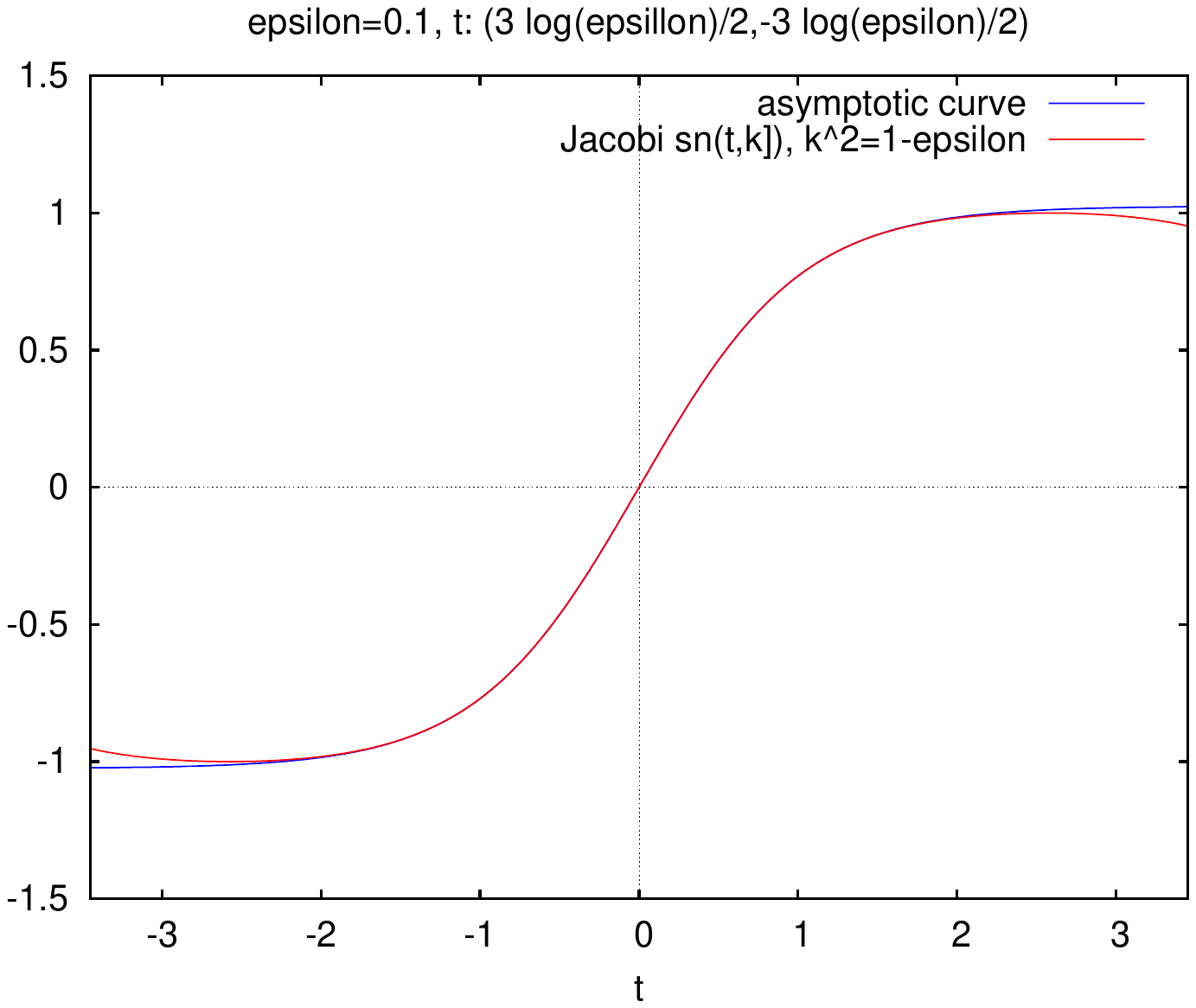}
\caption{Asymptotic curve and function $\sn(t|1-\epsilon)$ near the separatrix.}
\end{figure}

The constructed expansion is valid for less than a half of the period for the elliptic function $\sn$. This expansion is not valid for neighbourhoods of the turning points $u\sim1,u'=0$. 

To match this expansion to another which will be constructed  near the turning points we need an asymptotic properties of this expansion near the border of suitability.

Let us change the variable:

$$
t=-1/2\log(\epsilon)+\tau.
$$
Here $\tau$ is new independent variable. After substitution one obtains an asymptotic expansion as $\tau\ll-1$
\begin{eqnarray*}
u\sim 1+\epsilon\left(-\frac{e^{2\tau}}{128}-2 e^{-2\tau}+\frac{1}{4}\right)+
\\
\epsilon^2\left(\frac{\tau e^{2\tau}}{256}-\frac{\log(\epsilon)e^{2\tau}}{512} -\frac{5 e^{2\tau}}{512}-\tau e^{-2\tau} +\frac{\log(\epsilon)e^{-2\tau}}{2}-\frac{e^{-2\tau}}{2}+2e^{-4\tau}+\frac{11}{64}\right)
\end{eqnarray*}

The same asymptotic expansion can be obtained for neighbourhood of lower separatrix, if one uses the formula: $u(t+T/2,\epsilon)=-u(t,\epsilon)$. 

\subsection{Asymptotic behaviour near turning point} 

The elliptic function $\sn$ has an asymptotic behaviour of another type near the turning points. Here we constructs the asymptotic expansion for the $\sn$ near the  saddle-point $(1,0)$.
\begin{equation}
u(t,\epsilon)=1+\sum_{n=1}^\infty \epsilon^nv_n(\tau).
\label{asymptotoicsNearTurningPoint}
\end{equation}
Here $v_n=v_n(\tau)$.

The equations for the coefficients of the expansion can be obtained by ordinary way. One should collect  terms with the similar order of $\epsilon$. As a result one obtains a recurrent system of equations:
\begin{eqnarray*}
(v_1')^2=4v_1^2+2 v_1,
\\
2v_1'v_2'=8v_1 v_2-2v_2+4 v_1^3+5v_1,
\\
2v_1'v_n'=8v_1 v_n-2v_n+ P_n(v_1,\dots,v_{n-1}).
\end{eqnarray*}
Here $P_n$ is a polynomial of four power with $v_{k_1}v_{k_2} v_{k_3}v_{k_4}$, where $k_1+k_2+k_3+k_4=n$. 

\begin{figure}
\vspace{-9cm}
\includegraphics[scale=0.5]{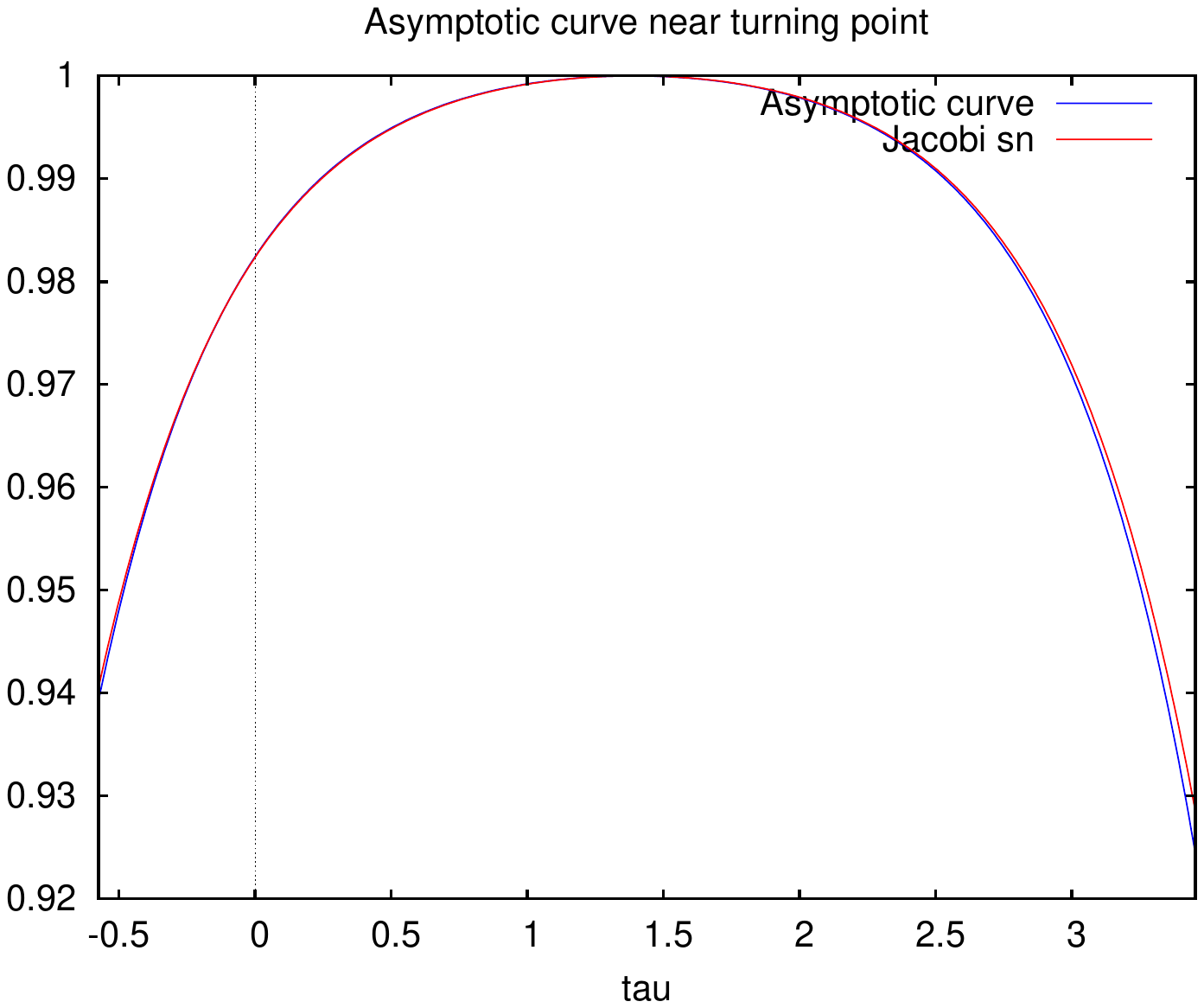}
\caption{The neighbourhood of the turning point $u=1$. The asymptotic curve and function $\sn(t-\log(\epsilon)/2,\sqrt{1-\epsilon}$ when $\epsilon=0.01$. The turning point not coincides with $\tau=0$, because for $\tau$ used only primary term of period for the elliptic function and not used the shift term  $2\log(2)$.
}
\end{figure}

The solution for $v_1$ has the form:  
$$
v_1=\frac{e^{2\tau}}{16} c_1+\frac{e^{-2\tau}}{4c_1}+\frac{1}{4}.
$$
Here $c_1$ is a parameter of solution.

The solution for the second-order term is:
$$
v_2=e^{2\tau}c_2-256 e^{-2\tau}c_2+\frac{e^{4\tau}}{32768}+\frac{\tau e^{2\tau}}{256}-\tau e^{-2\tau}-3 e^{-2\tau}+2 e^{-4\tau}+\frac{11}{64}.
$$
Here $c_2$ is a parameter of solution also.

The higher terms are solutions of linear equations of the first order. Their solutions can be presented in the form:
$$
v_n=e^{2\tau}c_n-256 e^{-2\tau}c_n+ O(e^{\pm2n\tau}),\quad \tau\to\pm\infty.
$$

Here $c_n$ is a parameter. It is defined by matching with asymptotic expansion which are valid outside of the small neighbourhoods of the turning points.

The validity of this expansion is defined by the condition:
$$
\epsilon^{n+1}v_{n+1}=o(\epsilon^{n} v_n).
$$
Using estimates for the grows of the terms one obtains:
$$
|\tau|\ll-1/2\log(\epsilon).
$$
The intervals of validity for (\ref{AsymptoticsCloseToSeparatrix}) are intersect when $t\gg1$ and $\tau\ll-1$. In the intersected field one can match the parameters of the asymptotic expansions. Here we choose the parameters $c_n$ of the terms $v_n$. In particular, the matching gives $c_1=-1/8$ and $c_2=-(\log(\epsilon)+5)/512$. 

As a result:
\begin{eqnarray*}
u(t,\epsilon)\sim 1-\epsilon\left(-\frac{e^{2\tau}}{128} -2e^{-2\tau}+\frac{1}{4}\right)+
\\
\epsilon^2
\left(
-e^{2\tau}\frac{(\log(\epsilon)+5)}{512}- e^{-2\tau}\frac{(\log(\epsilon)+5)}{2}+\right.
\\ 
\left. \frac{e^{4\tau}}{32768}+\frac{\tau e^{2\tau}}{256}-\tau e^{-2\tau}-3 e^{-2\tau}+2 e^{-4\tau}+\frac{11}{64} 
\right).
\end{eqnarray*}

\subsection{Uniform asymptotic expansion}

\begin{figure}[t]
\vspace{-8cm}
\includegraphics[scale=0.5]{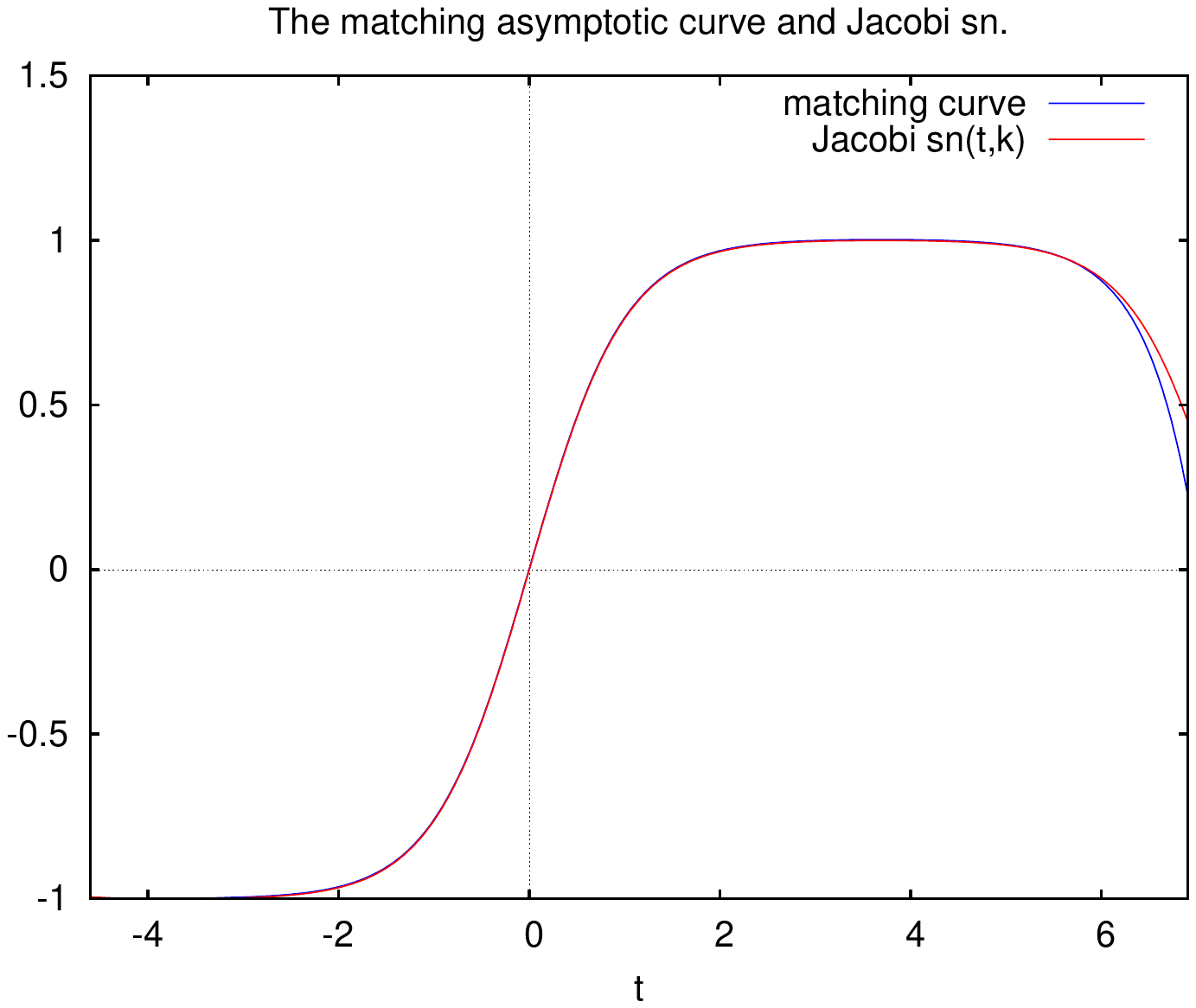}
\vspace{-1cm}
\caption{The combined asymptotic approximation, which is valid on over than a half of the period of elliptic function. $\epsilon=0.01$}
\label{fig-sn-matching-asymptotics-1}
\end{figure}

Now we are ready to construct a combined approximation of the function $u(t,\epsilon)\equiv\sn(t|1-\epsilon)$, which will be uniform over more than a half of the period  $t\in(\log(\sqrt{\epsilon}),-3\log(\sqrt{\epsilon}))$ as $\epsilon\to\infty$. To this we use an asymptotic device  which was offered by S.Kaplun \cite{Kaplun} for combined approximations.  We sum the constructed asymptotic expansions and subtract their common part. As a result we obtain the following formula (see figure \ref{fig-sn-matching-asymptotics-1}):
\begin{equation}
u(t,\epsilon)\sim \tanh(t)+\epsilon\frac{1}{\cosh^2(t)}\left(\frac{1}{8}\sinh(2t)-\frac{1}{4} t\right)+\frac{\epsilon}{4}-\epsilon^2\frac{e^{2t}}{128}.
\label{jacobi-sn-asymptotics-1}
\end{equation}

\begin{figure}[t]
\vspace{-8cm}
\includegraphics[scale=0.5]{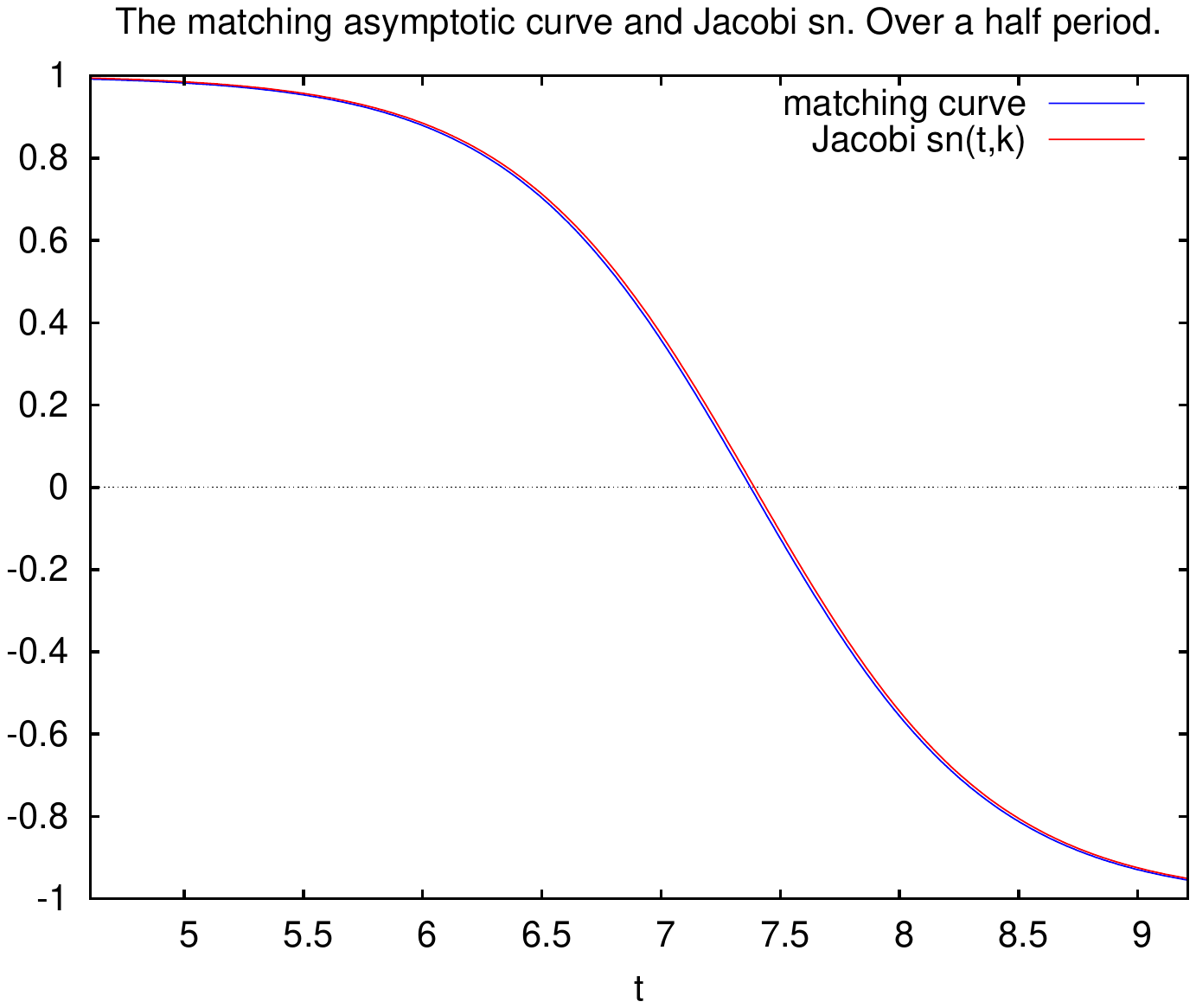}
\caption{
The combined asymptotic approximation which is valid on the second half of the period.$\epsilon=0.01$}
\label{fig-sn-matching-asymptotics-2}
\end{figure}

In the neighbourhood of the left saddle-point $u=-1$, $u'=0$ we can construct the same asymptotic expansion. But it is easy to get the formula  $u(t+T/2,\epsilon))=-u(t,\epsilon)$ and the asymptotic expansion will be obtained automatically way  using the asymptotic expansion (\ref{jacobi-sn-asymptotics-1}) (see figure  \ref{fig-sn-matching-asymptotics-2}).

\begin{eqnarray}
u(t,\epsilon)\sim -\tanh(t-T/2)-
\nonumber
\\
\epsilon\frac{1}{\cosh^2(t-T/2)}\left(\frac{1}{8}\sinh(2(t-T/2))-
\right.
\nonumber
\\
\left.
\frac{1}{4} (t-T/2)\right)+
\frac{\epsilon}{4}-\epsilon^2\frac{e^{2(t-T/2)}}{128}.
\label{jacobi-sn-asymptotics-2}
\end{eqnarray}

The combined asymptotic approximation which is valid over all period  of the elliptic function can be constructed by the same way, using the formulas (\ref{jacobi-sn-asymptotics-1}) and (\ref{jacobi-sn-asymptotics-2}) for the sum and subtract their common parts. But such combined asymptotic formula is large and does not written here.

\pagebreak

\end{document}